\documentclass[a4,12pt,multicol,makeidx]{article}
\def\origin{
  \clearpage
\vskip-\baselineskip\vskip-\topskip%
  \vbox to 0pt{\vskip-1in%
    \hbox to 0pt{\hskip-1in%
      \hbox to 0pt{\vrule width 1cm height .4pt depth 0mm\hss}%
      \vbox to 0pt{\hrule width .4pt height 0pt depth 1cm\vss}%
    \hss}%
  \vss}
  \vskip-\baselineskip
  \vbox to 0pt{\vskip-1in\vskip3cm%
    \hbox to 0pt{\hskip-1in\hskip3cm%
      \hbox to 0pt{\hss\vrule width 2cm height .4pt depth 0mm\hss}%
      \vbox to 0pt{\vss\hrule width .4pt height 1cm depth 1cm\vss}%
    \hss}%
  \vss}%
\vskip5mm\hskip10mm (3cm,3cm)
}%
   
\def\a{\alpha} 
\def\d{\delta} \def\l{\lambda}   \def\p{\partial}  \def\e{\varepsilon} 
\def\n{\nabla}    
\def\leq{\underline{<}} 
\def\la{\langle} \def\ra{\rangle} \def\hx{\hat{x}} 

\newenvironment{theorem}{%
\par \bigskip \it}{%
\bigskip \par}

\newenvironment{definition}{%
\par \bigskip \it}{%
\bigskip \par}

\makeindex
\title{A remark on the definitions of
viscosity solutions 
for the integro-differential equations with L{\'e}vy operators. 
}
\author{Mariko Arisawa\\ INRIA Rocquencourt
\\Domaine de Voluceau-Rocquencourt\\
BP 105, 78153 Le Chesnay, Cedex France\\
E-mail: mariko.arisawa@inria.fr
}
\date{}
\begin{document}
\maketitle
\bigskip
\section{Introduction.} 

$\quad$ In this note, we shall consider the following problem
\begin{equation}\label{first}
	F(x,u,\n u,\n^2 u)
	 -\int_{{\bf R^N}}  
	[u(x+z)
	-u(x)- {\bf 1}_{|z|\leq 1}\la \n u(x),z \ra]q(z) dz =0 \qquad x\in {\Omega},
\end{equation}
where $\Omega\subset {\bf R^N}$, $F$$\in C(\Omega\times {\bf R}\times {\bf R^N}\times {\bf S^N})$ is a second-order fully nonlinear elliptic operator, and the  L\'{e}vy measure $q(z)dz$ is a positive Radon measure such that 
\begin{equation}\label{integ}
	\int_{|z|< 1} |z|^2 q(z)dz+\int_{|z| \geq 1} 1 q(z)dz <\infty.
\end{equation}
The above type of problems is interested from the view point of 
the application in the mathematical finances (see Cont and Tankov \cite{cont}, Sulem and Oksendel \cite{os}). The comparison and the existence results have been studied in some frameworks of the viscosity solutions. However, the equivalence between these notions 
of viscosity solutions for (\ref{first}) are not trivial. Here, we would like to give some remarks on the relationships between 
 viscosity solutions defined in different manners. \\

	For an upper (resp. lower) semicontinuous function $u\in USC({\bf R^N})$ (resp.  $LSC({\bf R^N})$), we say that $(p,X)\in {\bf R^N}\times {\bf S^N}$ a subdifferential (resp. superdifferential)
 of $u$ at $x$, if for any $\d>0$ there exists $\e>0$ such that 
\begin{equation}\label{pX}
	u(x+z)-u(x)\leq    \quad \la p,z \ra + \frac{1}{2} \la Xz,z \ra + \d |z|^2 \quad \forall |z|\leq \e, 
\end{equation}
(resp.
\begin{equation}\label{pX2}
	u(x+z)-u(x) \geq  \quad \la p,z \ra + \frac{1}{2} \la Xz,z \ra - \d |z|^2 \quad \forall |z|\leq \e. 
\end{equation}
) 
We denote the set of all subdifferentials  (resp. superdifferentials) of $u$ at $x$  $J_{{\bf R^N}}^{2,+}u(x)$ (resp. $J_{{\bf R^N}}^{2,-}u(x)$).  
 As is well-known (see Crandall, Ishii and Lions \cite{users}), if $(p,X)$ is a subdifferential (resp. superdifferential) of $ u $ at $x$, then there exists $ \phi\in C^2(\mathbf{R^N}) $ 
 such that $ u(x)=\phi(x)$,  $ u-\phi$ takes a global maximum (resp. minimum) at $x$, and for any $\d>0$ there exists $\e>0$ such that 
\begin{equation}\label{phi}
	u(x+z)-u(x)\leq   
	\phi(x+z)-\phi(x)
	\leq    
	\la \n \phi(x),z \ra + \frac{1}{2} \la \n^2 \phi(x)z,z \ra + \d |z|^2 \quad \forall |z|\leq \e. 
\end{equation}
(resp.
\begin{equation}\label{phi2}
	u(x+z)-u(x)\geq    
	\phi(x+z)-\phi(x)
	\geq  
	\la \n \phi(x),z \ra + \frac{1}{2} \la \n^2 \phi(x)z,z \ra - \d |z|^2 \quad \forall |z|\leq \e. 
\end{equation}
)  
In Arisawa \cite{ar1}, \cite{ar2}, \cite{ar3}, the following definition of the viscosity solutions for (\ref{first}) was introduced. 
\begin{definition}{\bf Definition A.}  Let $u\in USC({\bf R^N})$ (resp. $v\in LSC({\bf R^N})$). We say that $u$ (resp. $v$) is a viscosity subsolution (resp. supersolution) of (\ref{first}), if for any $\hx\in \Omega$, any $(p,X)\in J_{{\bf R^N}}^{2,+}u(\hx)$ (resp. $\in J_{{\bf R^N}}^{2,-}v(\hx)$), and 
any pair of numbers $(\e,\delta)$ satisfying  (\ref{pX})  (resp. (\ref{pX2})), 
the following holds 
$$
	F(\hat{x},u(\hat{x}),p,X) 
	-\int_{|z|<\e} \frac{1}{2}\la(X+2\d I)z,z \ra q(z)dz\qquad\qquad\qquad\qquad\qquad\qquad
$$
\begin{equation}\label{def1}
	- \int_{|z|\geq \e} [u(\hat{x}+z)-u(\hat{x})
	-{\bf 1}_{|z|\leq 1} \la z,p\ra] q(z)dz \leq 0.
\end{equation}
(resp.
$$
	F(\hat{x},v(\hat{x}),p,X)
	-\int_{|z|<\e} \frac{1}{2}\la (X-2\d I)z,z\ra q(z)dz\qquad\qquad\qquad\qquad\qquad\qquad
$$
\begin{equation}\label{def2}	
	-\int_{|z|\geq \e} [v(\hat{x}+z)-v(\hat{x})
	-{\bf 1}_{|z|\leq 1} \la z,p\ra] q(z)dz  \geq 0.
\end{equation}
)  If $u$ is both a viscosity subsolution and a viscosity supersolution , it is called a viscosity solution.
\end{definition}
 We can rephrase Definition A by using the test functions in (\ref{phi}) (resp. (\ref{phi2})) as follows.

\begin{definition}{\bf Definition A'.}  Let $u\in USC({\bf R^N})$ (resp. $v\in LSC({\bf R^N})$). We say that $u$ (resp. $v$) is a viscosity subsolution (resp. supersolution) of (\ref{first}),  if for any $\hx\in \Omega$ and  for any $\phi\in C^2({\bf R^N})$ such that $u(\hx)=\phi(\hx)$ and $ u-\phi $ takes a global maximum (resp. minimum) at $\hx$, and for any pair of numbers $(\e,\d)$ satisfying (\ref{phi}) (resp. (\ref{phi2})), 
the following holds 
$$
	F(\hat{x},u(\hat{x}),\n \phi(\hat{x}),\n^2 \phi(\hat{x}))
	-\int_{|z|<\e} \frac{1}{2}\la(\n^2 \phi(\hat{x})+2\d I)z,z\ra q(z)dz
$$
\begin{equation}\label{def1'}
	- \int_{|z|\geq \e} [u(\hat{x}+z)-u(\hat{x})
	-{\bf 1}_{|z|\leq 1} \la z,\n \phi(\hat{x})\ra] q(z)dz \leq 0.
\end{equation}
(resp.
$$
	F(\hat{x},v(\hat{x}),\n \phi(\hat{x}),\n^2 \phi(\hat{x}))
	-\int_{|z|<\e} \frac{1}{2}\la ( \n^2 \phi(\hat{x}) -2\d I)z,z \ra q(z)dz
$$
\begin{equation}\label{def2'}	
	-\int_{|z|\geq \e} [v(\hat{x}+z)-v(\hat{x})
	-{\bf 1}_{|z|\leq 1} \la z,\n \phi(\hat{x})\ra] q(z)dz  \geq 0.
\end{equation}
)  If $u$ is both a viscosity subsolution and a viscosity supersolution , it is called a viscosity solution.
\end{definition}
 
We remark that the "global"  maximality (resp. minimality) of $u-\phi$ at $\hx$ in Definition A' can be replaced by the "local"  maximality (resp. minimality), 
 without changing any meaning of the definition. It is also clear that Definitions A and A' are equivalent. Next, we state the following definition of the viscosity solution  in Barles, Buckdahn and Pardoux \cite{bbp}, Jacobsen and Karlsen \cite{jk}, Barles and Imbert \cite{bi}.

\begin{definition}{\bf Definition B.}  Let $u\in USC({\bf R^N})$ (resp. $v\in LSC({\bf R^N})$). We say that $u$ (resp. $v$) is a viscosity subsolution (resp. supersolution) of (\ref{first}),  if for any $\hx\in \Omega$ and for any $\phi\in C^2({\bf R^N})$ such that $u(\hx)=\phi(\hx)$ and $ u-\phi $ takes a global maximum (resp. minimum) at $\hx$, 
\begin{equation}\label{defb1}
	F(\hat{x},u(\hat{x}),\n \phi(\hat{x}),\n^2 \phi(\hat{x}))
	- \int_{z\in \mathbf{R^N}} [\phi(\hat{x}+z)-\phi(\hat{x})
	- {\bf 1}_{|z|\leq 1} \la z,\n \phi(\hat{x})\ra] q(z)dz \leq 0.
\end{equation}
(resp.
\begin{equation}\label{defb2}	
	F(\hat{x},v(\hat{x}),\n \phi(\hat{x}),\n^2 \phi(\hat{x}))
	-\int_{z\in \mathbf{R^N} } [\phi(\hat{x}+z)-\phi(\hat{x})
	- {\bf 1}_{|z|\leq 1} \la z,\n \phi(\hat{x})\ra] q(z)dz \geq 0.
\end{equation}
)  
If $u$ is both a viscosity subsolution and a viscosity supersolution , it is called a viscosity solution.
\end{definition}

{\bf Remark 1.}  In the above cited works, Definition B was claimed to be equivalent to the following definition.

\begin{definition}{\bf Definition B'.}  Let $u\in USC({\bf R^N})$ (resp. $v\in LSC({\bf R^N})$). We say that $u$ (resp. $v$) is a viscosity subsolution (resp. supersolution) of (\ref{first}),  if for any $\hx\in \Omega$ and   for any $\phi\in C^2({\bf R^N})$ such that $u(\hx)=\phi(\hx)$ and $ u-\phi $ takes a global maximum (resp. minimum) at $\hx$, and for any $ \e>0$, 
$$
	F(\hat{x},u(\hat{x}),\n \phi(\hat{x}),\n^2 \phi(\hat{x}))
	-\int_{|z|<\e}   [\phi(\hat{x}+z)-\phi (\hat{x})
	-\la z,\n \phi(\hat{x})\ra]     q(z)dz
$$
\begin{equation}\label{defb1'}
	- \int_{|z|\geq \e} [u(\hat{x}+z)-u(\hat{x})
	-{\bf 1}_{|z|\leq 1}\la z,\n \phi(\hat{x})\ra] q(z)dz \leq 0.
\end{equation}
(resp.
$$
	F(\hat{x},v(\hat{x}),\n \phi(\hat{x}),\n^2 \phi(\hat{x}))
	-\int_{|z|<\e}  [\phi(\hat{x}+z)-\phi (\hat{x})
	-\la z,\n \phi(\hat{x})\ra]   q(z)dz
$$
\begin{equation}\label{defb2'}	
	-\int_{|z|\geq \e} [v(\hat{x}+z)-v(\hat{x})
	-{\bf 1}_{|z|\leq 1}\la z,\n \phi(\hat{x})\ra] q(z)dz \geq 0.
\end{equation}
)  
If $u$ is both a viscosity subsolution and a viscosity supersolution, it is called a viscosity solution.
\end{definition}
 The existence of the approximating sequence of test functions $ \phi_n(x) $ ($  \phi_n(x)\to  \phi (x)$ as $n\to \infty$, $ a.e.x$; 
$ u(x )\leq \phi_n(x)\leq  \phi(x)$ $\forall x\in \mathbf{R^N}$, $\forall n\in {\bf N}$,  in the case of the subsolution)  was used in the argument. 
Here, we shall consider Definition B, but not B'.\\

In this paper, thirdly we are interested in the following  definition of the viscosity solution, which seems to be  stronger than others at a first glance. 

\begin{definition}{\bf Definition C.}  Let $u\in USC({\bf R^N})$ (resp. $v\in LSC({\bf R^N})$). We say that $u$ (resp. $v$) is a viscosity subsolution (resp. supersolution) of (\ref{first}), if for any $\hx\in \Omega$ and  for any $\phi\in C^2({\bf R^N})$ such that $u(\hx)=\phi(\hx)$ and $ u-\phi $ takes a global maximum (resp. minimum) at $\hx$, the function $h(z)=u(\hat{x}+z)-u(\hat{x})-\la z,\n \phi((\hat{x})\ra$$\in L^1(\mathbf{R^N}, q(z)dz)$ and 
\begin{equation}\label{defc1}
	F(\hat{x},u(\hat{x}),\n \phi(\hat{x}),\n^2 \phi(\hat{x}))
	- \int_{z\in \mathbf{R^N}} [u(\hat{x}+z)-u(\hat{x})
	-{\bf 1}_{|z|\leq 1}\la z,\n \phi(\hat{x})\ra] q(z)dz \leq 0.
\end{equation}
(resp.
\begin{equation}\label{defc2}	
	F(\hat{x},v(\hat{x}),\n \phi(\hat{x}),\n^2 \phi(\hat{x}))
	-\int_{z\in \mathbf{R^N} } [v(\hat{x}+z)-v(\hat{x})
	-{\bf 1}_{|z|\leq 1}\la z,\n \phi(\hat{x})\ra] q(z)dz \geq 0.
\end{equation}
)  
If $u$ is both a viscosity subsolution and a viscosity supersolution , it is called a viscosity solution.
\end{definition}

We state the following results on the relationships between Definitions A, B and C.  \\

{\bf Theorem 1.$\quad$}
\begin{theorem} 
(i) If $u$ is the viscosity subsolution (resp. supersolution) of (\ref{first}) in the sense 
of Definition B, then $u$ is the viscosity subsolution (resp. supersolution) of (\ref{first}) in the sense 
of Definition C.\\
(ii) If $u$ is the viscosity subsolution (resp. supersolution) of (\ref{first}) in the sense 
of Definition C, then $u$ is the viscosity subsolution (resp. supersolution) of (\ref{first}) in the sense 
of Definition B.\\
\end{theorem}

{\bf Theorem 2.$\quad$}
\begin{theorem} 
(i) If $u$ is the viscosity subsolution (resp. supersolution) of (\ref{first}) in the sense 
of Definition A, then $u$ is the viscosity subsolution (resp. supersolution) of (\ref{first}) in the sense 
of Definition B.\\
(ii) If $u$ is the viscosity subsolution (resp. supersolution) of (\ref{first}) in the sense 
of Definition C, then $u$ is the viscosity subsolution (resp. supersolution) of (\ref{first}) in the sense 
of Definition A.\\
\end{theorem}

{\bf Theorem 3.$\quad$}
\begin{theorem} 
 The definitions A, B, and C are equivalent.
\end{theorem} 

In the following section 2, we first solve a technical problem, i.e. the construction of the sequence of test fuctions approximating 
 the subsolution $u$ from above. Then, in section 3 the above theorems will be proved by using the result of section 2. \\

We denote $B_{s}(x)=\{ y | \quad  | y-x | < s \}$$\subset  {\bf R^N}$ the ball centered at $x$ with the radius $s$, and 
$C_s (x)$$=\{ y| \quad  | y_i-x_i | < s \quad 1\leq i\leq N\}$$\subset  {\bf R^N}$ (where $x=(x_1,...,x_N)$, $y=(y_1,...,y_N)$) the cube centered at $x$ with the length of the edge $2s$. Moreover, we denote 
$$
	B_{s,s'}(x)= \{  y| \quad s< |y-x| <s' \}\subset  B_{s'}(x) \subset {\bf R^N}, 
$$
and denote 
$$
	C_{s,s'}(x)=\{ y| \quad s< | y_i-x_i | < s' \quad 1\leq i\leq N\}\subset  C_{s'}(x) \subset {\bf R^N}.
$$
 In the above notations, when $x=0$ we abbreviate as follows: $B_{s}=B_{s}(0)$,  $C_s=C_s (0)$, 
 $B_{s,s'}=B_{s,s'}(0)$, and $C_{s,s'}=C_{s,s'}(0)$.  
Let $P$ be a parallelotope which is the image of a linear transformation $T$ $(rank T=N)$ of a cube $C_s$, i.e. $ P=T C_s$. 
For $0<t<t' $ we denote 
$$
	P_t =T C_{ts}, \quad P_{t,t'} = TC_{ts,t's}. 
$$
We denote $P(x)=x+P$, $P_t(x)=x+P_t$, and  $P_{t,t'}(x)=x+P_{t,t'}$.

\section{Approximating sequence of test functions.}

$\quad$	Let $u(x)$ be an upper semi-continuous function. Assume that  there exists $\phi(x)\in C^2({\bf R^N})$, such that $u-\phi$ takes 
a global maximum at a point $\hx\in {\bf{R^N}}$ and $u(\hx)=\phi(\hx)$. In this situation,  we would like to construct  a sequence of test 
 functions $\psi_n\in C^2( {\bf{R^N}})$ $(n\in {\bf{N}})$, which 
roughly speaking, converges to $u$ as $n\to \infty$, by preserving 
 the following properties of $\phi$ at $\hx$: for any $n\in {\bf{N}}$\\
$$u-\psi_n\quad \hbox{takes a global maximum at }\quad \hx,$$
$$
	\n \psi_n(\hx)=\n \phi (\hx), \quad 
	\n^2 \psi_n(\hx)\geq \n^2 \phi (\hx), \quad \n^2 \psi_n(\hx)\downarrow \n^2 \phi (\hx)\quad \hbox{as}\quad n\to \infty.
$$

The construction of such a sequence is not trivial, and we obtained the following very near result. \\

{\bf Proposition 1.$\quad$} \begin{theorem} Let $u(x)\in USC(\bf R^N)$. Assume that there exists $\phi(x)\in C^2(\bf R^N)$, such that $u-\phi$ takes a global maximum at a point $\hx\in {\bf{R^N}}$ and $u(\hx)=\phi(\hx)$. Then, the following hold.\\

(i) For any $r\in (0,1)$ there exists $\psi^r(x)\in C^2({\bf R^N})$ and $P^r(\hx)$ a parallelotrope centered at  $\hx$, 
 such that 
$u-\psi^r$  takes a global maximum at $\hx$,
$u(\hx)=\psi^r(\hx)$, $\n \psi^r(\hx)=\n \phi (\hx$), $\n^2 \psi^r (\hx)= \n^2 \phi (\hx) + r I$, 
\begin{equation}\label{quadradef}
	\psi^r(x)=\psi^r(\hx)+ \la \n \phi(\hx), x-\hx \ra + \frac{1}{2}\la (\n^2 \phi(\hx) +rI )(x-\hx), x-\hx \ra
	\quad \hbox{in} \quad P^r_{\frac{1}{3}}(\hx),
\end{equation}
\begin{equation}\label{convex}
	\n^2 \psi^r(x)\geq O \quad \hbox{in}\quad  P^r_{\frac{2}{3},1}(\hx), 
\end{equation}
\begin{equation}\label{order}
	P^{r}(\hx) \subset P^{r'}(\hx) \quad \hbox{if} \quad r<r'; \quad   \lim_{r\to 0} \max_{x,y\in P_r(\hx)}|x-y|= 0.
\end{equation}
Moreover, there exists a constant $C>0$ independent on $r>0$ such that 
\begin{equation}\label{taylor}
	\psi^r(x)-\psi^r(\hx) -  \la \n \psi^r (\hx), x-\hx \ra \leq   C|x-\hx|^2.
	\quad \hbox{in}\quad P^r(\hx).
\end{equation}
(ii) For each $r>0$ there exists a sequence of functions $\psi^r_{n}(x)\in C^2({\bf R^N})$ ($n\in {\bf N}$) such that  
\begin{equation}\label{property}
	\psi^r_{n}(x) =\psi^{r}(x) \quad \hbox{in}\quad  P^r (\hx),  \quad 
	\lim_{n\to \infty} \psi^r_{n}(x) = u(x) \quad   \hbox{in}\quad  {\bf R^N}\backslash P^r(\hx), 
\end{equation}
\begin{equation}\label{monotones}
	\psi^r_{n+1}(x)<\psi^r_{n}(x) \quad \hbox{in}\quad {\bf R^N}\backslash P^r(\hx).
\end{equation}
\end{theorem}
$Proof$ $of$ $Proposition \quad 1.$ Without any loss of generality we may assume that $\hx=0$, $u(0)=\phi(0)=0$, $\n \phi(0)=0$. \\
(i) We shall first construct $\psi^r(x)$ for $r>0$. 
Put 
\begin{equation}\label{niji}
	\psi^r_{0} (x)= \la \n \phi(0),x \ra + \frac{1}{2} \la \n^2 \phi(0)x,x \ra+ \frac{r}{2}|x|^2.
\end{equation}
Since $\phi(0)=0$, there exists a number $s(r)>0$ such that 
$$
	\phi(x)\leq \psi^r_{0} (x) \quad x\in B_{2s(r)}, 
$$
 and that $u-\psi^r_{0}$ takes the global strict maximum at $0$ in $B_{2s(r)}$. We shall extend $\psi^r_{0}$ on ${\bf R^N}$, so that for the extended new function (by keeping the same notation) $\psi^r_{0}\in C^2({\bf R^N})$, $u-\psi^r_{0}$ takes its global strict maximum at $0$ in ${\bf R^N}$. 
Remark that the equation: $x_{N+1}=\psi^r_{0}(x)$ defined in $B_{2s(r)}\subset {\bf R^N}$ gives a quadratic surface in $B_{2s(r)}\times {\bf R}$$\subset {\bf R^{N+1}}$. Therefore, by the elementary 
result on the classification of the quadratic surface in the linear algebra,  by changing the coordinate system $x=(x_1,...,x_N)$ if necessary, 
 the quadratic surface given by (\ref{niji}) can be written in the following way
\begin{equation}\label{eigen}
	\psi^r_{0} (x)=\sum_{i=1}^{N} \l_i x_i^2 \quad x\in B_{2s(r)},
\end{equation}
where $\l_i$ ($1\leq i\leq N$) are the eigenvalues of the matrix $\n^2 \phi(0)+ rI$, we still use the notation $(x_1,...,x_N)$ for the new coordinate system, 
 and the equation may be considered to hold in $B_{2s(r)}$ of the new coordinate system.  
We need the following lemma. \\

{\bf Lemma 1.$\quad$} \begin{theorem}
Let $\l < 0$, $s>0$, and consider $f(x)=\l x^2$ in the interval $-s\leq x \leq s$. Let 
 $g(x)=a \exp(-\frac{c}{|x-\a|^2})+b$, where $\a=\frac{2s}{3}$, $a=-\frac{e \l  s^2}{9}$, $b=\frac{2\l s^2}{9}$, and $c=\frac{s^2}{9}$.  Define 
$$
	\psi(x)= f(x) \quad 0\leq x\leq \frac{s}{3};\quad =g(x) \quad  \frac{s}{3}\leq x < \frac{2s}{3}; \quad =\frac{2\l s^2}{9} \quad  \frac{2s}{3}\leq x\leq s, 
$$
and $\psi(x)=\psi(-x)$ ($-s\leq x \leq 0$). Then, $\psi(x)$ is $C^2$ in $-s\leq x \leq s$,  $\psi(x)= f(x)$ in $|x|\leq  \frac{s}{3}$ and  
$\psi(x)$ is convex in  $\frac{2s}{3}\leq |x|\leq  s$. 
\end{theorem}
$Proof$ $of$ $Lemma \quad 1.$ By the elementary calculation, we see that $f(\frac{s}{3})=g(\frac{s}{3})$$=\frac{\l s^2}{9}$,  
 $f'(\frac{s}{3})=g'(\frac{s}{3})$$=\frac{2\l s}{3}$, $f''(\frac{s}{3})=g''(\frac{s}{3})$$=2\l$, $\lim_{x\uparrow \frac{2s}{3}}g(x)=\frac{2\l s^2}{9}$, 
$\lim_{x\uparrow  \frac{2s}{3}}g'(x)=\lim_{x\uparrow   \frac{2s}{3}}g''(x)=0$. 
Thus, we get the function $\psi$ as in the claim. \\

Assume that $\l_i<0$ ($1\leq i \leq n$), and $\l_i\geq 0$ ($n+1 \leq i \leq N$) in (\ref{eigen}). 
Remark that $C_{s(r)}\subset B_{2s(r)}$, and 
by using the above lemma for $\l=\l_i$ ($1\leq i \leq n$), put
 $$
	\psi_i(x)=\psi(x_i) \quad \hbox{for}\quad x\in C_{s(r)}\subset {\bf R^N}.
$$  Define 
\begin{equation}\label{base}
	\psi^r(x)=\sum_{i=1}^{n}\psi_i(x) + \sum_{i=n+1}^{N}\l_i x_i^2 \quad  \hbox{for}\quad x\in C_{s(r)}.
\end{equation}
Then from Lemma 1, $\psi^r(x)=\psi^r_{0}(x)$ in $C_{\frac{s(r)}{3}}$, $\psi^r(x)$ is convex in $C_{\frac{2s(r)}{3},s(r)}$. Consider now the original coordinate 
system, by putting $P^r=TC_{s(r)}$, where $T$ represents the linear transformation to the original coordinate system.  The above argument leads (\ref{quadradef}) and (\ref{convex}) in the corresponding parallerotope  $P^r_{\frac{1}{3}}$, and the doughnut type region $P^r_{\frac{2}{3},1}$.  As for (\ref{order}), if $r<r'$ holds then we can take $s(r)<s(r')$, and the claim is clear from the above argument.  
From Lemma 1, 
$$
	0 \geq \psi_i(x) \geq \min_{|x_i|< s(r)} \l_i x_i^2 \quad x\in C_{s(r)}, \quad 1\leq i \leq n. 
$$ 
Therefore, 
by taking account the way that $\psi^r(x)$ is constructed from $\psi^r_{0}(x)$ (quadratic in $C_{s(r)}$) in (\ref{eigen}) and (\ref{base}), it is clear that the following holds. 
$$
	\psi^r(x)-\psi^r(0) -  \left\langle \n \psi^r(0) ,x \right\rangle  \leq max_{1\leq i\leq N} \l_i |x|^2 \leq C|x|^2 \quad x\in C_{s(r)},
$$
where $C=|\n^2 \phi(0)|+1$. 
We consider the above inequality in the original coordinate system, and  see that (\ref{taylor})  holds in  $P^r= T C_{s(r)}$. 
 Therefore, we get the function  $\psi^r(x)$ in (i).\\

(ii) Let $\rho_n>0$ ($n\in {\bf N}$) be a sequence of numbers such that $\lim_{n\to \infty} \rho_n=0$.
 From (i) $\psi^r(x)$ is convex in $P^r_{\frac{2}{3},1}$, and 
$\psi^r(x) > u(x)$ in $P^r\backslash \{0\}$. Thus, for each $n\in {\bf N}$ we can extend $\psi^r(x)$ on $P^r_{1+\rho_n}$$(\supset P^r)$ so that the extended function 
 $\psi^r_{n}(x)$ is $C^2$,  satisfying 
$$
	\psi^r_{n}(x)=\psi^r(x) \quad\hbox{in}\quad P^r, \quad \psi^r_{n}\quad \hbox{is convex in}\quad P^r_{\frac{2}{3},1+\rho_n}, 
$$
 and $u-\psi^r_{n}$ takes a global strict maximum at $0$ in $P^r_{1+\rho_n}$. Furthermore, since $P^r_{1+\rho_{n+1}}\subset P^r_{1+\rho_n}$, we can 
 extend $\psi^r_{n}$ on ${\bf R^N}$ ($n\in {\bf N}$) so that the extended functions (keeping the same notations) $\psi^r_{n}(x)\in C^2({\bf R^N})$, and 
\begin{equation}\label{a}
	\psi^r_{n+1}(x)<\psi^r_{n}(x) \quad \hbox{in}\quad {\bf R^N}\backslash P^r, \quad \forall n\in {\bf N}, 
\end{equation}
\begin{equation}\label{b}
	\lim_{n\to \infty} \psi^r_{n}(x) = u(x) \quad \hbox{in}\quad  {\bf R^N}\backslash P^r. 
\end{equation}
 Remark that  (\ref{a}),  (\ref{b}) are possible, for $\psi^r_{n}$ ($n\in {\bf N}$) are convex on $\p P^r$. 
Therefore, we have constructed the sequence  $\psi^r_n (x)$ ($n\in {\bf N}$)  in (ii). \\
 
If we do not need the convergence of  the second-order derivatives of the test functions: 
$\n^2 \psi_n(\hx)\downarrow \n^2 \phi (\hx)\quad \hbox{as}\quad n\to \infty$, the construction of the approximating sequence 
 is much simpler. The idea of the following comes from a result in Evans \cite{evans}.\\

{\bf Proposition 2.$\quad$} \begin{theorem} Let $u(x)\in USC({\bf R^N})$. Assume that there exists  $\phi(x)\in C^2(\bf R^N)$, such that $u-\phi$ takes a global maximum at a point $\hx\in {\bf{R^N}}$ and $u(\hx)=\phi(\hx)$.  
 Then, there exists a sequence of functions $\psi_n(x)$$\in C^2(\mathbf{R^N})$  such that $u-\psi_n$ 
takes the global maximum at $\hat{x}$, $u(\hat{x})=\psi_n(\hat{x})$,  $\n \phi(\hat{x})=\n \psi_n(\hat{x})$, 
 and 
$$
	\lim_{n\to \infty} \psi_n ({x})=u(x),\quad \psi_n (x)\geq  u(x) \quad \forall x\in \mathbf{R^N}.
$$
\end{theorem}
$Proof$ $of$ $Proposition \quad 2.$  We may assume that $\hat{x}=0$, $u(\hx)=\phi(\hx)=0$, $\n \phi(\hat{x})=0$, without any loss of the generality.  
Now, since $\phi\in C^2$, we can take $M_n=\sup_{|x|\leq n^{-1}} |\n^2 \phi (0)|$ for any $n\in \mathbf{N}$. Put 
$\psi^0_n (x)=2M_n |x|^2$ in $\{|x|\leq n^{-1}\}$, and extend it to $\mathbf{R^N}$ so that $\psi^0_n (x)\geq \phi(x)$, $\psi^0_n (x)\in C^2$ on $\mathbf{R^N}$. Remark that $ \psi^0_n-u$ takes its global maximum at $0$ for any $n\in {\bf N}$. Since $ \psi^0_n$ is convex and radially symmetric in $\{|x|\leq n^{-1}\}$, 
 we can take $ \psi_n$ such that 
$$
	\psi_n(x)= \psi^0_n(x) \quad \hbox{for} \quad |x|\leq (2n)^{-1};  \quad 
	u(x)\leq \psi_n(x) \leq u(x)+n^{-1}\quad \hbox{for} \quad |x|\geq 2n^{-1},  
$$
$$
	 \psi_{n+1}(x) \leq  \psi_n(x)  \quad \hbox{on} \quad \mathbf{R^N}, \quad \hbox{for} \quad \forall n\in {\bf N}.
$$
The sequence of functions $\{\psi_n\}$ ($n\in {\bf R^N}$) satisfies the claim, clearly. \\

{\bf Remark 2.} From the above construction of $\psi_n(x)$, we only have 
$$
	\n^2 \phi(\hat{x})\leq \n^2 \psi_n(\hat{x}) \quad  \hbox{for} \quad  \forall \quad  n.
$$
 Proposition 2 can be used to prove the equivalence of the definitions of viscosity solutions for (\ref{first}), when $F$ is the first-order Hamiltonian.\\

{\bf Remark 3.} The construction of the approximating sequence of test functions for the supersolution can be done similarly.\\

\section{Proofs of the main results.}

$\quad$We use the following well-known elementary theorem of the monotone convergence of Beppo-Levi.\\

{\bf Lemma$\quad$2.} (Beppo-Levi, see H. Brezis \cite{bre}.) 
\begin{theorem} Let $f_n(x)$ ($n\in {\bf N}$) be a sequence of increasing functions in $L^1(\mathcal{O}, d\mu (x))$ ($\mathcal{O}\in {\bf R^N}$), such that 
 $\sup_n \int_O f_n d\mu (x) < \infty$. Then, $f_n(x)$ converges almost everywhere in $\mathcal{O}$ to a function $f(x)$. Moreover 
$f(x)\in L^1$ and $|| f_n - f ||_{L^1}\to 0$ as< $n\to \infty$. 
\end{theorem}

We begin with the proof of Theorem 1.\\

$Proof$ $of$ $Theorem$ $1.$ (i) Let $u$ be a viscosity subsolution (resp. supersolution) of (\ref{first}) in the sense 
of Definition B. Assume that there exists $\phi\in C^2(\mathbf{R^N})$ such that 
$u-\phi$ takes a global maximum at ${\hx}\in \Omega$,  and 
$u({\hx})=\phi({\hx})$. Let $r>0$ be an arbitrary small number. Then from Lemma 1,  there exists a parallelotrope $P^r(\hx)$, a function $\psi^r\in C^2$, and a 
sequence of functions $\psi^r_{n}\in C^2$ ($n\in {\bf N}$) having the properties in (i), (ii) of Lemma 1. Since 
$u-\psi^r_{n}$ ($n\in {\bf N}$) takes a global maximum at $\hx$, 
 from Definition B
$$
	F({\hx},u({\hx}),\n \psi^r_{n}({\hx}), \n^2 \psi^r_{n}({\hx}))\qquad\qquad\qquad\qquad\qquad\qquad\qquad
$$
$$
	- \int_{z\in \mathbf{R^N}} [\psi^r_{n}({\hx}+z)-\psi^r_{n}({\hx})
	-{\bf 1}_{|z|\leq 1}\la z,\n \psi^r_{n}({\hx})\ra] q(z)dz  \leq 0 \quad \forall n .
$$
From Lemma 1 (ii) (\ref{property}), the above can be written as follows.
$$
	F({\hx},u({\hx}),\n \psi^r_{n}(\hx), \n^2 \psi^r_{n}({\hx}))\qquad\qquad\qquad\qquad\qquad\qquad\qquad
$$
$$
	- \int_{\hx + z \in {P^r(\hx)}} [\psi^{r}({\hx}+z)-\psi^{r}({\hx})
	-{\bf 1}_{|z|\leq 1}\la z,\n \psi^{r}({\hx})\ra] q(z)dz  
$$
$$
	- \int_{\hx + z \in (P^r(\hx))^c} [\psi^r_{n}({\hx}+z)-\psi^r_{n}({\hx})
	-{\bf 1}_{|z|\leq 1} \la z,\n \psi^r_{n}({\hx})\ra] q(z)dz  \leq 0 \quad \forall n.
$$
Put 
$$
	h_n(z)=\psi^r_{n}({\hx}+z)-\psi^r_{n}({\hx})
	-{\bf 1}_{|z|\leq 1}\la z,\n \psi^r_{n}({\hx})\ra \quad \forall n. 
$$ 
Then, from the continuity of $F$ and (\ref{taylor}), we have 
$$
	\sup_n [- \int_{  \hx + z \in (P^r(\hx))^c }  h_n(z) q(z)dz]
$$
$$
	\leq \sup_n [ -F( {\hx},u({\hx}),\n \psi^r_{n}(\hx), \n^2 \psi^r_{n}({\hx}) ) 
	+ C \int_{\hx + z \in {P^r(\hx)}} |z|^2 q(z) dz] < \infty. 
$$
From (\ref{property}), (\ref{monotones}), $h_n(z)$ is monotone decreasing as $n\to \infty$ and 
$$
	\lim_{n\to \infty} h_n(z) = u({\hx}+z)-u({\hx})-{\bf 1}_{|z|\leq 1}\la z,\n \phi({\hx})\ra \quad z\in \{z|  \hx + z \in (P^r(\hx))^c \}. 
$$
Thus, from Lemma 2 (Beppo-Levi) we see $u({\hx}+z)-u({\hx})-{\bf 1}_{|z|\leq 1}\la z,\n \phi({\hx})\ra$$\in L^1({\bf R^N},q(z)dz)$, and 
$$
	- \int_{  \hx + z \in (P^r(\hx))^c } [u({\hx}+z)-u({\hx})
	-{\bf 1}_{|z|\leq 1}\la z,\n \phi({\hx})\ra] q(z)dz  \quad\qquad\qquad\qquad
$$
$$
	\leq 
	-F({\hx},u({\hx}),\n \phi({\hx}), \n^2 \phi({\hx})+rI)  
	+ C_r < \infty, 
$$
where $C_r>0$ is a constant such that $\lim_{r\to 0}C_r=0$ from (\ref{order}). 
Now, by letting $r\to 0$ in the above inequality, from the continuity of $F$ and (\ref{order}) 
$$
	F({\hx},u({\hx}),\n \phi({\hx}), \n^2 \phi({\hx}))
	- \int_{z\in \mathbf{R^N}} [u({\hx}+z)-u({\hx})
	-{\bf 1}_{|z|\leq 1} \la z,\n \phi({\hx})\ra] q(z)dz \leq 0 
$$
holds. 
Therefore, $u$ is the viscosity subsolution in the sense of Definition C.\\

(ii) Let $u$ be a viscosity subsolution (resp. supersolution) of (\ref{first}) in the sense 
of Definition C.  Assume that there exists $\phi\in C^2(\mathbf{R^N})$ such that 
$u-\phi$ takes a global maximum at ${\hx}\in \Omega$,  and 
$u({\hx})=\phi({\hx})$. From Definition C, 
$$
	F({\hx},u({\hx}),\n \phi({\hx}),\n^2 \phi({\hx}))
	- \int_{z\in \mathbf{R^N}} [u({\hx}+z)-u({\hx})
	-{\bf 1}_{|z|\leq 1} \la z,\n \phi(({\hx})\ra] q(z)dz  \leq 0.
$$
Since $u(\hx+z)\leq \phi(\hx+z)$ for any $z\in \mathbf{R^N} $, it is clear hat the above leads 
$$
	F({\hx},u({\hx}),\n \phi({\hx}),\n^2 \phi({\hx}))
	- \int_{z\in \mathbf{R^N}} [\phi({\hx}+z)-\phi({\hx})
	-{\bf 1}_{|z|\leq 1} \la z,\n \phi({\hx})\ra] q(z)dz \leq 0.
$$
Therefore,  $u$ is the viscosity subsolution in the sense of Definition B.\\

{\bf Remark 4.} If $F$ is the first-order Hamiltonian, the approximating sequence $\psi_n$ ($n\in {\bf N}$) in Proposition 2 serves 
 to prove the claim  in Theorem 1. \\

Next, we shall prove Theorems 2 and 3.\\

$Proof$ $of$ $Theorem \quad 2.$ (i) Let $u$ be a viscosity subsolution (resp. supersolution) of (\ref{first}) in the sense 
of Definition A. Remark that Definition A is equivalent to Definition A'. Assume that there exists $\phi\in C^2(\mathbf{R^N})$ such that 
$u-\phi$ takes a global maximum at ${\hx}\in \Omega$,  and 
$u({\hx})=\phi({\hx})$. Then, for any pair of numbers $(\e,\d)$ such that (\ref{phi}) holds, 
$$
	F({\hx},u({\hx}),\n \phi({\hx}),\n^2 \phi({\hx}))
	-\int_{|z|<\e} \frac{1}{2}\la(\n^2 \phi({\hx})+2\d I)z,z\ra q(z)dz
$$
$$
	- \int_{|z|\geq \e} [u({\hx}+z)-u({\hx})
	-{\bf 1}_{|z|\leq 1}\la z,\n \phi({\hx})\ra] q(z)dz \leq 0.
$$
Then, since $u(\hx+z)\leq \phi(\hx+z)$ for any $z\in \mathbf{R^N} $, 
$$
	F({\hx},u({\hx}),\n \phi({\hx}),\n^2 \phi({\hx}))
	-\int_{|z|<\e} \frac{1}{2}\la(\n^2 \phi({\hx})+2\d I)z,z\ra q(z)dz
$$
$$
	- \int_{|z|\geq \e} [\phi ({\hx}+z)-\phi ({\hx})
	-{\bf 1}_{|z|\leq 1}\la z,\n \phi({\hx})\ra] q(z)dz\leq 0.
$$
By tending $\e \to 0$, this shows that $u$ is the viscosity solution in the sense of Definition B. \\

(ii) Let $u$ be a viscosity subsolution (resp. supersolution) of (\ref{first}) in the sense 
of Definition C. Assume that there exists $\phi\in C^2(\mathbf{R^N})$ such that 
$u-\phi$ takes a global maximum at ${\hx}\in \Omega$,  and 
$u({\hx})=\phi({\hx})$. We have 
$$
	F({\hx},u({\hx}),\n \phi({\hx}),\n^2 \phi({\hx}))
	- \int_{z\in \mathbf{R^N}} [u({\hx}+z)-u({\hx})
	-{\bf 1}_{|z|\leq 1}\la z,\n \phi({\hx})\ra] q(z)dz  \leq 0.
$$
Since 
$$
	u(\hx+z)-u(\hx)-\la z,\n \phi({\hx})\ra \leq \phi(\hx+z)-\phi(\hx)-\la z,\n \phi({\hx})\ra
$$
$$
	\leq 
	\frac{1}{2} \la \n^2 \phi(\hx)z,z \ra + \d |z|^2 \quad |z|\leq \e, 
$$
we have 
$$
	F({\hx},u({\hx}),\n \phi({\hx}),\n^2 \phi({\hx}))
	-\int_{|z|<\e} \frac{1}{2}\la(\n^2 \phi({\hx})+2\d I)z,z\ra q(z)dz
$$
$$
	- \int_{|z|\geq \e} [u({\hx}+z)-u({\hx})
	-{\bf 1}_{|z|\leq 1} \la z,\n \phi({\hx})\ra] q(z)dz  \leq 0.
$$
That is, Definition  C implies Definition A. \\

{\bf Remark 5.} For the viscosity supersolutions, the similar claims to those in Theorems 1 and 2 hold, too.\\

$Proof$ $of$ $Theorem \quad 3.$ The claim comes directly from Theorems 1 and 2.\\


\end{document}